# Integrating reliability and resilience to support microgrid design

*Raphael Wu[1], Giovanni Sansavini[1,*]*

[1]*Reliability and Risk Engineering Laboratory, Institute of Energy Technology, Mechanical and Process Engineering Department, ETH Zurich, Switzerland*

* Corresponding Author: e-mail: sansavig@ethz.ch

**Abstract**

Quantifying the potential benefits of microgrids in the design phase can support the transition of passive distribution networks into microgrids. At current, reliability and resilience are the main drivers for this transition. Therefore, this paper presents a mathematical optimization model to support the retrofitting of distribution networks into microgrids integrating techno-economic, resilience and reliability objectives. Storage and distributed generation are optionally installed to complement renewable generation, enabling the microgrid to supply priority demands during stochastic islanding events with uncertain duration. For a comprehensive quantification and optimization of microgrid resilience and reliability, islanding due to external events is combined with a detailed model of internal faults. Minimizing the interruption costs yields optimal capacities and placements of distributed energy resources and new lines for reconfiguration. The proposed method produces microgrid designs with up to 95% reliability and resilience gain and moderate cost increase in two benchmark distribution networks using data from the US Department of Energy. The developed methodology is scalable to large networks owing to the tailored Column-and-Constraint-Generation approach.

**Nomenclature**

| | |
|---|---|
| $E_i$ | CCG subproblem constraint matrix for master variables |
| $G_i$ | CCG subproblem constraint matrix for subproblem variables |
| $h_i$ | CCG subproblem constraint right hand sides |
| $I$ | Set of all islanding subproblems |
| $I_{inf}$ | Set of infeasible islanding subproblems |
| $I_m$ | Set of subproblems added to the master problem |
| $LB$ | Lower bound to the optimal CCG objective |
| $O_i$ | Optimal objective function value of subproblem i |
| $\mathbb{P}_{PCC,Dmd}$ | Path between PCC and $Dmd$ bus |

| $t$ | time index |
|---|---|
| $T$ | Time horizon length |
| $UB$ | Upper bound to the optimal CCG objective |
| $X_i$ | Feasible region for CCG subproblem variables |
| $X_m$ | Feasible region for CCG master problem variables |

**Acronyms**

| | |
|---|---|
| CCG | Column-and-Constraint-Generation |
| EENS | Expected Energy Not Supplied |
| GEV | Generalized Extreme Value distribution |
| MED | Major Event Day |
| MG | Microgrid |
| MILP | Mixed Integer Linear Program |
| PV | Photovoltaic |
| SAIDI | System Average Interruption Duration Index |
| SAIFI | System Average Interruption Frequency Index |

**Sub- and Superscripts**

| | |
|---|---|
| b | Bus |
| ch | Storage charging |
| d | Storage discharging |
| DER = {DG,RG,Str} | Distributed Energy Resources |
| DG | Distributed Generator |
| Dmd | Demand bus |
| ex | Export |
| i | Islanding event |
| im | Import |
| Inv | Investment |
| l | Line |
| m | Master problem |
| Op | Operation |
| PCC | Point of Common Coupling |
| Rel | Reliability |
| Res | Resilience |
| RG | Renewable generator |
| self | Storage self-discharge |
| Str | Storage |

**Parameters**

| | |
|---|---|
| $\alpha$ | Yearly interest rate |
| $a_l, a^{DER}$ | Annuity factor for lines and DER |
| $c_{Inv,0}^{DER}$ | Fixed DER investment cost |
| $c_{Inv,S}^{DER}$ | Variable DER investment price |
| $c_P^{DER}$ | DER active power price |
| $c_Q^{DER}$ | DER reactive power price |
| $c_{NS}^{Dmd}$ | Cost of energy not supplied at bus Dmd |
| $c_{l,Inv}$ | Line investment cost |
| $c_m$ | Cost coefficients in CCG master problem |
| $c_{ex,t}^{PCC}$ | PCC power export price |
| $c_{im,t}^{PCC}$ | PCC power import price |
| $c_Q^{PCC}$ | PCC reactive power price |
| $c_{Curt}^{RG}$ | RG curtailment price |
| $d_i$ | CCG subproblem objective coefficients |
| $DOD_{max}$ | Storage maximum depth-of-discharge |
| $E_i^{Dmd}$ | Energy demand in islanding event i |
| $E_i^{RG}$ | Renewable energy available in islanding event i |
| $\eta_{self}, \eta_{ch}, \eta_d$ | Storage self-discharge-, charge- and discharge efficiency |
| $f_{P,max}^{DG}$ | DG active power limit |
| $f_{Q,max}^{DG}$ | DG reactive power limit |
| $f_C^{Str}$ | Ratio of storage power to energy rating |
| $\lambda^{Dmd}$ | Failure rate of equipment at bus Dmd |
| $\lambda_l$ | Line failure rate |
| $M$ | Large number for "Big-M" constraints |
| $N_0$ | Number of islanding events considered in the first CCG iteration |
| $N(b)$ | Set of lines adjacent to bus b |
| $n_{Cyc,max}$ | Maximum allowed number of storage charge/discharge cycles |
| $n_p$ | Number of polygon sides in a polygonal inner approximation |
| $P_t^{Dmd}$ | Active power demand |
| $P_{t,Avl}^{RG}$ | Available renewable power |
| $p_i$ | Islanding event occurrence probability |
| $p_{t_i}$ | probability of islanding event i to end at time $t_i$ |
| $PF_{min}^{DG}$ | DG minimum power factor |
| $PF_t^{Dmd}$ | demand power factor |
| $Q_t^{Dmd}$ | Reactive power demand |
| $r_l$ | line resistance |
| $S_t^{Dmd}$ | Apparent power demand |
| $S_{l,max}$ | Line rated power |
| $t_i$ | Time index of islanding event i |
| $t_{Lf}$ | DER or line life time |
| $\tau^{Dmd}$ | Repair time of equipment at bus Dmd |
| $\tau_l$ | Line repair duration |
| $\tau_{max}$ | Maximum repair duration of all lines |
| $V_{b,min}, V_{b,max}$ | Minimum and maximum squared bus voltage magnitude |
| $w_t$ | Time step weight |
| | |
| **Optimization variables** | |
| $C_{Inv}$ | Investment cost |
| $C_{Op,t}$ | Operations cost |
| $C_{Res,t}$ | Resilience cost |
| $C_{Rel}$ | Reliability cost |
| $C_{Op,t}^{DER}$ | DER operations cost |
| $C_{Op,t}^{PCC}$ | PCC operations cost |
| $C_{Curt,t}^{RG}$ | RG curtailment cost |
| $C_{i,t_i}^{Str}$ | Cost of recharging storage after islanding event i |
| $EENS^{Dmd}$ | Expected energy not supplied at bus Dmd |
| $E_{max}^{Str}$ | Storage energy capacity |
| $E_t^{Str}$ | Stored energy |
| $\eta$ | Worst-case islanding event objective |
| $P_t^{DER}$ | DER active power |
| $P_{net,t}^{Dmd}$ | Net demand after subtracting local supply at bus Dmd |
| $P_t^{DG}$ | DG active power |
| $P_{l,t}$ | Line active power |
| $P_t^{PCC}$ | Net PCC active power import including losses |
| $P_{DistFl,t}^{PCC}$ | Net PCC active power import without losses |
| $P_t^{RG}$ | RG active power |
| $P_t^{Str}$ | Storage net discharging power |
| $P_{ch,t}^{Str}, P_{d,t}^{Str}$ | Storage charging and discharging power |
| $Q_t^{DER}$ | DER reactive power |
| $Q_t^{DG}$ | DG reactive power |
| $Q_{l,t}$ | Line reactive power |
| $Q_{DistFl,t}^{PCC}$ | PCC reactive power import |
| $Q_{Pos,t}^{PCC}$ | Absolute value of PCC reactive power |
| $Q_t^{Str}$ | Storage reactive power |
| $S_{max}^{DER}$ | DER power capacity |
| $S_{max}^{DG}$ | DG power capacity |
| $S_{max}^{RG}$ | RG power capacity |
| $S_{max}^{Str}$ | Storage power capacity |
| $U_{Rel}^{Dmd}$ | Expected outage duration at bus Dmd |
| $V_{b,t}$ | Squared bus voltage magnitude |
| $x_i$ | Variables of islanding event i |
| $x_m$ | CCG master problem variables |
| $x_m^*$ | Optimal value of $x_m$ at the current CCG iteration |
| $y_{Inv}^{DER}$ | Binary DER investment variable |
| $y_l^{Dmd}, \bar{y}_l^{Dmd}$ | Binary variables modelling the path between Dmd and the PCC |
| $y_l$ | Binary line investment variable |

## 1. Introduction

Microgrids (MGs) have the capability to operate in islanding mode owing to distributed energy resources (DER), thus potentially ensuring reliability and resilience in electricity distribution [1]. Furthermore, controllable DER can balance fluctuating renewable generation (RG) and increase RG hosting capacity [2]. Financial benefits compared to passive grids arise from saving grid use fees, peak shaving, reducing losses, and deferring grid upgrades [3]. However, due to the high investment costs of DER, the aforementioned benefits are fully unlocked if the design of MGs takes the transition from the current mostly passive power distribution system into due account [2]. To this aim, this paper investigates the optimal transformation of distribution networks with RG into reliable and resilient MGs.

*1.1. Literature review*

This literature review focuses on optimal design of MGs including resilience and reliability objectives, which identifies a complex engineering task and an active field of research. General reviews of MGs and their optimal design are provided in [2], [4], [5].

*Resilience and reliability in distribution grids*

Electricity distribution grids and microgrids are called resilient if they can withstand extreme events [6]–[8]. In [6], islanding events are defined as power outages lasting longer than a specified threshold of 5 hours and up to 24 hours, whereas shorter outages are categorized as reliability events. Hussain et al. [7] refer to resilience as the ability to avoid load shedding during extreme weather events and review effective operation strategies to enhance resilience. They identify proactive operation strategies to prepare for extreme events considering the probability of occurrence and uncertain duration, and the consideration of storage state of charge as research gaps. These are addressed in the present paper.

Jufri et al. [8] distinguish between reliability and resilience: reliability is defined as the ability to withstand common outages, which are events with high probability, short duration and small impact regions. Reliability study is conducted using outage frequency and durations, and quantified with the indices System Average Interruption Frequency Index (SAIFI) and System Average Interruption Duration Index (SAIDI). Resilience assessment, on the other hand, should include extreme event anticipation, absorption and recovery, with an active change in grid operation to mitigate event impacts.

This paper adopts the definitions in [8], where resilience is the ability to anticipate, absorb and recover from extreme events with long durations, whereas reliability is the ability to supply demands despite common outages of short durations.

*Islanding in MG design*

Extended islanding operations are the unique feature of MGs and must be included in their design [1]. Islanding requirements are often represented as adequacy constraints, ensuring that the total generation capacity exceeds demand power in the MG at every time [9], [10]. This implementation, however, may lead to overestimating system resilience because voltage and line flow limits may require load shedding during islanding [11], [12].

Furthermore, storage, though vital for MGs to integrate RG and provide backup power during islanding [13], is often excluded from MG design optimization [10], [12] or it is modeled via exogenous charging cycles [9]. This avoids analyzing islanding events, during which storage operation should differ from grid-connected mode. Ref. [11] is an exception, where storage is optimally placed and operated, considering islanding during the most costly 8 hours of the year. Although the storage can supply all demands during islanding, operation costs are underestimated because storage reserves are not maintained at all times to cope with the uncertain timing of islanding events [7]. This leads to increased costs of distributed generators (DG) or grid imports and produces suboptimal solutions. The short islanding durations of 8 hours or less [11], [14] also indicates that the designed MGs are not necessarily resilient to long-lasting disruptions, despite their ability to islanding. Therefore, optimizing storage design and operation for islanding events with uncertain timing remains an open challenge.

Most MGs use chemical storage technologies [2], [5], [6], [13], [15], [16]. Lead-acid batteries are a mature technology with 70-80% efficiency [17]. Redox flow systems are highly efficient and can scale power and energy ratings independently, but are still expensive [17]. Lithium based batteries are highly efficient with decreasing cost and increasing cycle life [13]. Sodium-sulfur batteries which operate around 350°C achieve high efficiency and reasonable cycle life [17]. Flywheels and capacitors can supply high power for short periods of time [17]. Compressed air and pumped hydro are rarely used due to local resource requirements [18]. This paper employs a storage model which can be parameterized to represent a large variety of generic storage types [16].

*Reliable MG design*

Even though faults within MGs can influence the optimal design [14], they are rarely considered. Switches and DER are co-optimized to reduce internal fault impacts in [14]. In [19], a loop topology is imposed to have redundant supply paths for each demand, but without DER optimization. Although [15] combines islanding and faults within the MG by scaling demand bus failure rates by the probability of MGs having to import power, the resulting MG designs are not resilient according to the definitions of [6]–[8]. The absence of an islanding event model prevents the MG from anticipating or changing its dispatch while islanding, and the short event durations and high frequencies are within the realm of reliability rather than resilience [8]. Therefore, the optimal design of reliable and resilient MGs is still an open research gap.

*Optimization methods*

MG design literature can be further classified by the optimization method used. Mathematical optimization used in [11], [12], [19] results in a solution with optimality guarantee [20], but requires approximations of the power flow equations [16]. To solve computationally demanding operating scenarios, Benders decomposition is used to separate the investment and operation level in [11]. Heuristic [10] and metaheuristic [9], [14], [15], [21] approaches can use the exact power flow formulation but provide no optimality guarantee [20].

## *1.2. Contributions and paper organization*

This paper develops a new model to investigate the optimal design of MGs for economic, resilience and reliability objectives by combining the siting and sizing of DER and power lines with operational optimization. The paper addresses the following research questions:

1. How can islanding events with uncertain timing and duration be incorporated into optimal MG design?
2. How to transition distribution grids with RG optimally to reliable and resilient MGs?
3. What is the influence of reliability, resilience and their combination on optimal MG design?
4. What is the sensitivity of optimal MG design to different islanding events?

The paper's main contributions are:
- modeling resilience against stochastic islanding events with uncertain timing and duration, including DER redispatch and power flow;

- proposing a reliability model for faults within the MG, where DER can supply their local bus, considering energy, active and apparent power limitations;
- combining resilience, reliability, design and grid-tied operation in a single optimization framework, enabling the optimal design of MGs which are able to anticipate and mitigate the effects of islanding events [8];
- adapting a Column-and-Constraint-Generation method (CCG) [22] to solve the design problem with significant computational improvements and scalability to large systems with respect to stochastic optimization.

This paper is organized as follows: Section 2 introduces the model and solution method. The method is applied to two case studies described in Section 3, the results of which are reported in Section 4. Section 5 draws conclusions and answers the research questions.

## 2. Microgrid model and optimization methodology

Section 2.1 presents an overview of the optimization, followed by a full description in Sections 2.2 – 2.6. Section 2.7 describes the adapted CCG method used to solve the optimization.

### 2.1. Optimization Framework and objective

The proposed MG design optimization problem (1) optimizes the design and operations of MGs, considering resilience and reliability:

$$\min \left[ C_{Inv} + \sum_{t=1}^{T} w_t C_{Op,t} + \sum_{t=1}^{T} w_t C_{Res,t_{i,1}=t} + C_{Rel} \right] \quad (1)$$
$$subject\ to\ (2) - (39).$$

The objective function is the expected equivalent annual cost [23] of the MG and consists of the investment costs $C_{Inv}$, the operations costs $C_{Op,t}$, the resilience costs of islanding events $C_{Res,t_{i,1}=t}$ and the reliability costs from faults within the MG $C_{Rel}$. To reduce the computational effort, the representative day method is used [24], with $w_t$ as the number of represented days in the year.

To facilitate the model presentation, the optimization problem (1) is structured into four levels, displayed in Fig. 1. Despite the hierarchical structure, the whole optimization (1) is solved as a single problem, considering all levels simultaneously.

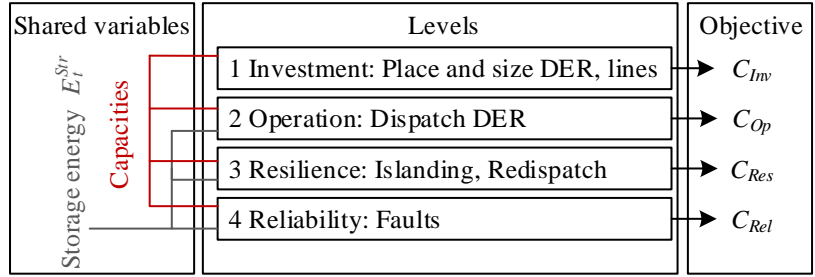

*Fig. 1. Optimization framework overview. Connections on the left-hand side indicate decision variables that are shared between levels.*

Decision variables, objective and constraints are attributed to investment (Level 1, Section 2.2), operations (Level 2, Section 2.3), resilience (Level 3, Section 2.4) and reliability (Level 4, Section 2.5). Within the optimization problem (1), the levels are linked by shared decision variables pertaining to investment in DER and in power lines, and by the storage energy profile $E_t^{Str}$ as shown on the left side of Fig. 1.

## 2.2. Investment

The investment level contains variables and constraints related to the location and capacity of power lines and DER.

$$C_{Inv} = \sum_{DER} \frac{c_{Inv,S}^{DER} S_{max}^{DER} + c_{Inv,0}^{DER} y_{Inv}^{DER}}{a^{DER}} + \sum_l \frac{c_{l,Inv} y_l}{a_l} \quad (2)$$

$$0 \leq S_{max}^{DER} \leq M y_{Inv}^{DER} \ \forall \ DER \quad (3)$$

$$y_{Inv}^{DER}, y_l \in \{0,1\} \ \forall \ DER, l. \quad (4)$$

The DER cost in Eq. (2) consists of a capacity-dependent part $c_{Inv,S}^{DER}$ and a fixed part $c_{Inv,0}^{DER}$, multiplied by the installed capacity $S_{max}^{DER}$ and a binary variable $y_{Inv}^{DER}$ representing the installation status, respectively. Line costs are determined using fixed costs $c_{l,Inv}$ for each candidate line, and a binary variable $y_l$ indicating that line $l$ has been installed. All investment costs are converted to equivalent annual costs using the annuity factor $a = \frac{1-(1+\alpha)^{-t_{Lf}}}{\alpha}$, where $t_{Lf}$ is the life time of a DER or line, and $\alpha$ is the yearly interest rate [23]. Setting a large $M$ in Eq. (3) ensures that as soon as a nonzero DER capacity is installed, the installation status $y_{Inv}^{DER}$ is equal to 1 [25].

### 2.3. Operations

*Distributed Generators*

The DG model represents inverter-based and directly coupled DG, which can provide active ($P_t^{DG}$) and reactive power ($Q_t^{DG}$) [26]:

$$(P_t^{DG})^2 + (Q_t^{DG})^2 \leq (S_{max}^{DG})^2 \ \forall \ t, DG \quad (5)$$
$$0 \leq P_t^{DG} \leq f_{P,max}^{DG} S_{max}^{DG} \ \forall \ t, DG \quad (6)$$
$$-f_{Q,max}^{DG} * S_{max}^{DG} \leq Q_t^{DG} \leq f_{Q,max}^{DG} S_{max}^{DG} \ \forall \ t, DG \quad (7)$$
$$|Q_t^{DG}| \leq P_t^{DG} * \tan \operatorname{acos} PF_{min}^{DG} \ \forall, DG \quad (8)$$

Apparent power magnitudes are limited by inverter or stator limits, depending on the DG type (5). Furthermore, there can be active (6), reactive (7) limits where $f_{P,max}^{DG}$ and $f_{Q,max}^{DG} \in [0, 1]$ are DG type-dependent factors. A minimum power factor $PF_{min}^{DG}$, defined as the minimum ratio of active and apparent power magnitudes $|P_t^{DG}|/|S_t^{DG}|$, is used to enforce technical or regulatory limits [27] in (8).

Renewable generators (RG) such as PV cells and wind turbines are a special case of DG and also include an active power limit, where the available active power $P_{t,Avl}^{RG}$ is determined by weather conditions and conversion efficiencies, taking PV module efficiency and DC-AC losses into account:

$$P_t^{RG} \leq P_{t,Avl}^{RG} \ \forall \ t, RG. \quad (9)$$

*Battery Storage*

The battery storage (Str) model presented in (10) – (17) improves the model in [16] by including self-discharge losses :

$$(P_t^{Str})^2 + (Q_t^{Str})^2 \leq (S_{max}^{Str})^2 \ \forall \ t, Str \tag{10}$$

$$S_{max}^{Str} = f_C^{Str} E_{max}^{Str} \ \forall \ Str \tag{11}$$

$$(1 - DOD_{max}) E_{max}^{Str} \leq E_t^{Str} \leq E_{max}^{Str} \ \forall t, Str \tag{12}$$

$$E_t^{Str} = \eta_{self} E_{t-1}^{Str} + \eta_{ch} P_{ch,t}^{Str} - \frac{P_{d,t}^{Str}}{\eta_d} \ \forall \ t > 1, Str \tag{13}$$

$$P_t^{Str} = P_{d,t}^{Str} - P_{ch,t}^{Str} \ \forall t, Str \tag{14}$$

$$P_{d,t}^{Str} \geq 0, P_{ch,t}^{Str} \geq 0 \ \ \forall t, Str \tag{15}$$

$$E_{t=1}^{Str} = E_{t=24*k}^{Str} \ \forall \ Str, k = [1..T/24] \tag{16}$$

$$\sum_t P_{d,t}^{Str} + P_{ch,t}^{Str} \leq 2n_{Cyc,max} E_{max}^{Str} \ \forall \ Str \ . \tag{17}$$

Analogous to (5) for DG, (10) limits storage active and reactive power to the inverter's rated capacity, which is proportional to the battery-cell energy capacity (11). The useful storage energy is limited in (12) by the capacity and the maximum depth of discharge $DOD_{max}$. The energy balance (13) considers self-discharge ($\eta_{self}$) and charge/discharge ($\eta_{ch}, \eta_d$) losses. $P_t^{Str}$ is positive when discharging (14). Eq. (15) combined with the objective function prevents simultaneous charging and discharging, which would waste energy. Eq. (16) ensures equal storage energy levels at the beginning and end of each day with hourly time steps. Aging is considered by limiting the number of charge/discharge cycles to $n_{Cyc,max}$ (17); this avoids excessive degradation during the expected life time. With no loss of generality, it is assumed that PV is installed by building owners and battery storage is installed by the MG designer. Therefore, PV and storages have separate inverters. Similar to [16], there is no dedicated battery control algorithm. Battery charging and discharging is optimized through the overall cost function (1).

*Power flow and power exchange through PCC*

Distribution grids are commonly operated radially, therefore, the linearized DistFlow equations are used [16], modeling active ($P_{l,t}$) and reactive ($Q_{l,t}$) power flows as well as squared voltage magnitudes $V_{b,t}$, which are incorporated into power balance (18) – (19), voltage (20) and line flow limits (21):

$$P_t^{Str} + P_t^{DG} + P_t^{RG} - P_t^{Dmd} = \sum_{l \in N(b)} P_{l,t} \; \forall \, t, b \quad (18)$$
$$Q_t^{Str} + Q_t^{DG} + Q_t^{RG} - Q_t^{Dmd} = \sum_{l \in N(b)} Q_{l,t} \; \forall \, t, b \quad (19)$$
$$V_{b,min} \leq V_{b,t} \leq V_{b,max} \; \forall \, t, b \quad (20)$$
$$(P_{l,t})^2 + (Q_{l,t})^2 \leq y_l (S_{l,max})^2 \; \forall \, t, l \quad (21)$$

Eqs. (18) – (19) show the most general case, where all DER technologies and demands are present at bus $b$.

In contrast to [16], resistive losses cannot be multiplied by the electricity cost and added to the objective function, as the MG could be entirely self-supplied or exporting energy. Therefore, the active power import at the PCC is modeled as:

$$P_t^{PCC} = P_{DistFl,t}^{PCC} + \sum_l r_l (P_{l,t}^2 + Q_{l,t}^2) \; \forall \, t \quad (22)$$

where $P_{DistFl,t}^{PCC}$ is the active power import to the MG computed by the lossless DistFlow equations [16] and the second term approximates the losses with line resistances $r_l$. In the operational objective, $P_t^{PCC}$ is charged with electricity costs in the objective function. Therefore, the MG can compensate losses internally by decreasing $P_{DistFl,t}^{PCC}$.

Finally, the contribution of the operations to the objective function in Eq. (1) consists of costs for electricity from the PCC and DER, as well as renewable curtailment fees:

$$C_{Op,t} = C_{Op,t}^{PCC} + \sum_{DER} C_{Op,t}^{DER} + \sum_{RG} C_{Curt,t}^{RG} \; \forall \, t \quad (23)$$
$$C_{Op,t}^{DER} = (c_P^{DER} P_t^{DER} + c_Q^{DER} |Q_t^{DER}|) \Delta t \; \forall \, t, DER \quad (24)$$
$$C_{Curt,t}^{RG} = c_{Curt}^{RG} (P_{t,Avl}^{RG} - P_t^{RG}) \Delta t \; \forall \, t, RG \quad (25)$$
$$C_{Op,t}^{PCC} = (c_{ex,t}^{PCC} P_t^{PCC} + (c_{im,t}^{PCC} - c_{ex,t}^{PCC})(P_t^{PCC} + |P_t^{PCC}|)/2 \\ + c_Q^{PCC} |Q_{DistFl,t}^{PCC}|) \Delta t \; \forall \, t \quad (26)$$

DER operations and maintenance costs in (24) include costs for active power production $c_P^{DER}$ (e.g. fuel costs for DG), and the production or consumption of reactive power $c_Q^{DER}$ (e.g. inverter losses [27]). For customer-owned renewable generators, curtailment is compensated with an additional fee $c_{Curt}^{RG}$ in (25) [27].

Prices for electricity consumption are often higher than for feed-in. Eq. (26) ensures that the price $c_{ex,t}^{PCC}$ is used for active power export and $c_{im,t}^{PCC}$ for active power import, assuming $c_{im,t}^{PCC} \geq c_{ex,t}^{PCC} \ \forall \ t$. Using $P_t^{PCC}$ instead of $P_{DistFl,t}^{PCC}$ in (26) enables the MG to compensate losses internally, by decreasing $P_{DistFl,t}^{PCC}$ to $-\sum_l r_l (P_{l,t}^2 + Q_{l,t}^2)$.

For reactive power, an equal price $c_Q^{PCC}$ is assumed for import and export, and the reactive power $Q_{DistFl,t}^{PCC}$ computed by the lossless DistFlow equations is used to calculate the absolute value $|Q_{DistFl,t}^{PCC}|$.

### 2.4. Resilience

Following the resilience concepts in [6]–[8], islanding events are modeled with a stochastic approach including uncertain starting time and duration, as shown in Fig. 2.

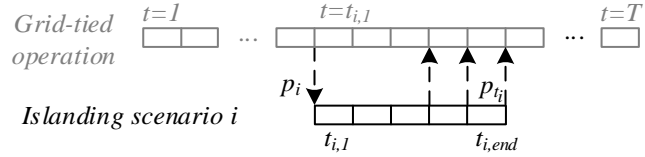

Fig. 2. Visualization of a stochastic islanding event.

Islanding event $i$ starts at time $t_{i,1}$ with probability $p_i$ and its uncertain duration is represented by $p_{t_i}$, i.e. the probability of the islanding event ending at any time $t_i = [t_{i,1}, ..., t_{i,end}]$. The maximum islanding duration $t_{i,end}$ is given by the time step at which the cumulative sum of ending probabilities is equal to 1. As MGs should be able to island during unforeseen disturbances [7], a separate islanding event is considered for every time step $t$ of the optimization.

Every islanding event $i$, starting at time $t_{i,1} = t$ is modeled with a separate set of operational variables and constraints (5) – (25) to enable a redispatch when the MG is islanded, combined with the following additional constraints:

$$P_{i,t_i}^{PCC} = Q_{DistFl,i,t_i}^{PCC} = 0 \ \forall \ i, t_i \qquad (27)$$

$$E_{i,t_{i,1}}^{Str} = E_{t=t_{i,1}}^{Str} \ \forall \ i, t_i, Str \qquad (28)$$
$$C_{i,t_i}^{Str} = c_{im,t=t_{i+1}}^{PCC}\left(E_{t=t_i}^{Str} - E_{i,t_i}^{Str}\right) \forall \ i, Str \qquad (29)$$
$$C_{Res,i} = p_i \sum_{t_i} p_{t_i} \left(\sum_{Str} C_{i,t_i}^{Str} + \sum_{\tau=t_{i,1}}^{t_i} \sum_{DER} C_{Op,i,\tau}^{DER}\right) \forall \ i. \qquad (30)$$

During islanding, the power exchange at the PCC is zero (27). Storage systems are initialized at the energy content they have at the beginning of the contingency (28), i.e. time $t_{i,1}$. During islanding, the storage energy profile $E_{i,t_i}^{Str}$ can differ from its value $E_{t=t_i}^{Str}$ for grid-connected operations. After reconnecting to the grid at the end of an islanding event, the difference between the two storage levels is charged with the grid import fee $c_{im,t=t_{i+1}}^{PCC}$ in (29). This assumes that the stored energy is brought back to the grid-tied operational optimum level after reconnection to the grid using imports from the PCC. The storage recovery costs $C_{i,t_i}^{Str}$ following the islanding event are combined with the costs of DER redispatch $C_{Op,i,\tau}^{DER}$ in the islanding objective $C_{Res,i}$ in (30), which quantifies the expected value of DER redispatch and recovery costs. As the proposed method aims at designing MGs which can withstand islanding events without load shedding, all demands must be supplied while islanding. Embedding this islanding model in the overall framework introduced in Section 2.1 ensures that the MG is operated in anticipation of islanding events which can happen at any time step $t$, considering the probability of event occurrence $p_i$ and the uncertain duration characterized by $p_{t_i}$. Embracing the resilience definition as the ability to anticipate, absorb and recover from extreme events with an active change in grid operation [6]–[8], modeling islanding stochastic islanding events with uncertain starting time and duration make the MG resilient because: 1) the MG prepares for extreme events causing a blackout by being able to island at any time, 2) if the event occurs, the MG changes its operation by ramping up DER to supply all loads, mitigating the impact; 3) recovery is modeled with the cost of buying back depleted storage energy (29).

## 2.5. Reliability

To include MG reliability against common outages characterized by failure rates and durations [8], faults within the MG are modeled based on the approach in [28]: Every line $l$ is characterized by a failure rate $\lambda_l$ and a repair duration $\tau_l$. Additional equipment such as cable splices and transformers installed at demand bus $Dmd$ is characterized by $\lambda^{Dmd}$ and $\tau^{Dmd}$. The expected duration $U_{Rel}^{Dmd}$ of outages at every demand bus $Dmd$ is determined by [28]:

$$\sum_{l \in N(b)} y_l^{Dmd} - \bar{y}_l^{Dmd} = \begin{cases} -1; b = Dmd \\ 1; b = PCC \\ 0; otherwise \end{cases} \forall b, Dmd \quad (31)$$

$$y_l^{Dmd}, \bar{y}_l^{Dmd} \in \{0,1\} \ \forall \ l, Dmd \quad (32)$$

$$U_{Rel}^{Dmd} = \lambda^{Dmd}\tau^{Dmd} + \sum_l (y_l^{Dmd} + \bar{y}_l^{Dmd})\lambda_l \tau_l \ \forall \ Dmd. \quad (33)$$

Only lines $l$ belonging to the path $\mathbb{P}_{PCC,Dmd}$ from the supply (PCC) and $Dmd$ are contributing to the outage frequencies and durations at bus $Dmd$. In case of multiple paths, the optimization chooses the path resulting in the lowest failure rate when minimizing unreliability. To identify whether line $l$ belongs to $\mathbb{P}_{PCC,Dmd}$, two sets of binary variables $y_{l,Dmd}$ and $\bar{y}_{l,Dmd}$ are introduced and constrained using (31) – (32), where $N(b)$ is the set of lines adjacent to bus $b$. The sum $y_l^{Dmd} + \bar{y}_l^{Dmd}$ takes a value of 1 if line $l$ belongs to $\mathbb{P}_{PCC,Dmd}$, and zero, otherwise [28]. Therefore, $U_{Rel}^{Dmd}$ can be expressed with (33).

To compute demand not supplied due to faults, the expected outage duration and mean demand are needed [29]. With appropriate grid codes and control, DER such as storage [26], RG and DG [30] can self-supply their bus, even if the rest of the grid is not energized. Therefore, the approach in [28] is extended by including self-supply for DER host buses during upstream faults. $P_{net,t}^{Dmd}$ models the net demand after subtracting local supply, which is limited by the installed capacity (34) and the available active power (35) – (36) at each bus:

$$P_{net,t}^{Dmd} \geq PF_t^{Dmd}\left(\left|S_t^{Dmd}\right| - S_{max}^{RG} - S_{max}^{DG} - S_{max}^{Str}\right) \forall\, t, Dmd \tag{34}$$

$$P_{net,t}^{Dmd} \geq P_t^{Dmd} - P_{t,Avl}^{RG} - f_{P,max}^{DG} S_{max}^{DG} - E_t^{Str} \eta_d / \sum_{k=1}^{\tau_{max}} \eta_{self}^{1-k} \forall\, t, Dmd \tag{35}$$

$$P_{net,t}^{Dmd} \geq P_t^{Dmd} - P_{t,Avl}^{RG} - f_{P,max}^{DG} S_{max}^{DG} - S_{max}^{Str} \forall\, t, Dmd \tag{36}$$

$$P_{net,t}^{Dmd} \geq 0 \; \forall\, t, Dmd. \tag{37}$$

Eqs. (34) – (36) show the most general case, where all types of DER are installed at bus $Dmd$. Apparent (34) and active (35) – (36) power limits are considered separately due to varying RG availability and storage energy. As storage active power is limited by the energy content $E_t^{Str}$ (35) and the inverter capacity $S_{max}^{Str}$ (36), two constraints model the limits on active power. To get a conservative estimate of the available storage energy at time $t$, the term $E_t^{Str} \eta_d / \sum_{k=1}^{\tau_{max}} \eta_{self}^{1-k}$ is the available active power in 1 time step if the storage is emptied at constant power over the next $\tau_{max}$ time steps, where $\tau_{max} = \max_l \tau_l$.

The expected energy not supplied (EENS) at each demand bus is calculated as:

$$EENS^{Dmd} = U_{Rel}^{Dmd} \sum_t w_t P_{net,t}^{Dmd} / T \; \forall\, Dmd \tag{38}$$

$$C_{Rel} = \sum_{Dmd} c_{NS}^{Dmd} EENS^{Dmd} \tag{39}$$

The reliability objective contribution $C_{Rel}$ in (39) is the cost of EENS in the entire MG, with $c_{NS}^{Dmd}$ as cost of demand not supplied.

### 2.6. Linearization

Four types of nonlinearities have to be replaced to convert (1) to a MILP, which can be efficiently solved to global optimality. Quadratic constraints (5), (10) and (21) are replaced by a polygonal inner approximation [31]. The square terms in (22) are approximated with a piecewise linear approximation [32]. Products of binary and continuous variables in (38) are replaced by an exact linear reformulation [32]. Absolute values $X^+$ of variables $X$ in (8), (24), and (26) are modeled as:

$$X^+ \geq X \text{ and } X^+ \geq -X \tag{40}$$

Although (40) does not guarantee tightness ($X^+ = |X|$) by itself, the linearization is only used if $X^+$ is in the objective function with a positive coefficient, or when limiting $X$. In both cases, (40) is either tight at the optimum, or irrelevant as $X$ has not reached the limits imposed by other constraints.

### 2.7. Optimization solution method

Due to the detailed reliability and resilience models, solving (1) is computationally extremely challenging for any realistic case study. Therefore, a dedicated solution method must be developed. The linearized optimization problem (1) is a stochastic MILP with one islanding event per time step $t$ and is solved with a modified CCG method. To provide a clear introduction, the modified CCG method is presented first. Subsequently, the improvements with respect to the original CCG method [22] are highlighted and motivated.

The modified CCG minimizes the worst-case cost of a problem under uncertainty using the process shown in Fig. 3, which encompasses the following iterative steps:

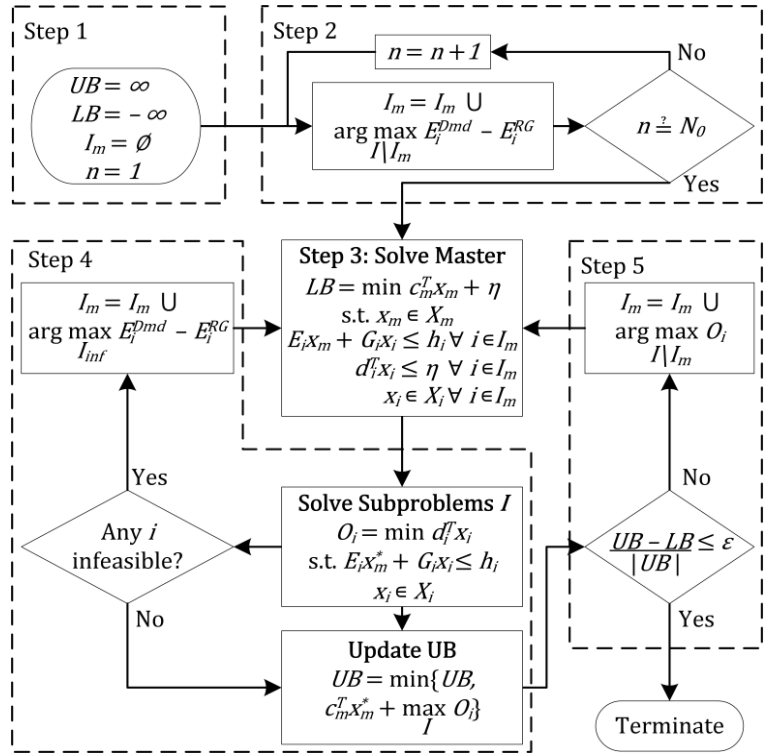

*Fig. 3. Column-and-Constraint-Generation procedure*

1. the problem is split into a master problem with variables $x_m$, constraints $x_m \in X_m$ and objective $c_m^T x_m$, and one subproblem for each islanding event $i \in I$ (Level 3 in Fig. 1), with variables $x_i$, constraints $E_i x_m + G_i x_i \leq h_i, x_i \in X_i$ and objective $d_i^T x_i$. Following the notation in Section 2.4, the subproblem objective is $C_{Res,i}$, nonzero rows in $E_i$ are the capacity and storage energy continuity constraints (28) – (29), and $G_i$ contains the islanding operation constraints. The upper ($UB = \infty$) and lower bound ($LB = -\infty$) on the objective, and the set $I_m = \emptyset$ of subproblems considered when solving the master problem are initialized;

2. $N_0 \geq 1$ subproblems are added to $I_m$, prioritizing the subproblems with maximum demand $E_i^{Dmd}$ after subtracting previously installed RG energy $E_i^{RG}$;

3. the master problem including all events in $I_m$ is minimized. An objective term $\eta$ is added and constrained to be greater than or equal to all subproblems' objectives within $I_m$. As the master problem only considers a subset of islanding events, its optimized objective value provides a lower bound to the original problem's worst-case objective;

4. given the current optimal solution of the master problem $x_m^*$ computed in Step 3, all subproblems are solved separately. If any subproblem is infeasible, the subproblem with maximum $E_i^{Dmd} - E_i^{RG}$ from the infeasible set of subproblems $I_{inf}$ is added to $I_m$, and Step 3 is repeated. If all problems are feasible, an upper bound is calculated by summing the master objective $c_m^T x_m^*$ and the maximum subproblem objective $\max_I O_i$;

5. if the relative difference between upper and lower bound is below a predefined threshold $\varepsilon$, the algorithm has converged. Otherwise, the subproblem $i$ with the largest objective function value is added to $I_m$ and the procedure returns to Step 3.

This modified CCG method improves the computational efficiency of the original approach [22] in two ways:

First, islanding events with maximum net demand are added in Steps 2 and 4, which speeds up convergence because the islanding event cost and the necessary DER capacity to provide the energy are correlated with the net demand $E_i^{Dmd} - E_i^{RG}$. Thus, the probability of having infeasible or more costly subproblems in subsequent iterations is decreased. In the original CCG method [22], subproblems are added in an arbitrary order.

Second, the original CCG method solves the master problem without any subproblems in the first iteration and defines $N_0 = 0$. When solving (1) with $N_0 = 0$, the first few CCG iterations are spent adding islanding events during the peak demand day, as they impose capacity requirements to enable islanding. With $N_0 > 1$, DER are installed and several high demand islanding events are feasible after the first iteration.

As the highest demand islanding events are typically adjacent in time during the peak demand day, defining $N_0 > 1$ also serves as a non-anticipation constraint for the energy storage; indeed, if $N_0 = 1$, the optimization considers only the peak demand event in the first iteration, filling the storage just before the beginning of this event. With $N_0 > 1$, the storage is operated to enable islanding for all of the $N_0$ highest demand events.

3. **Case studies**

The developed method is applied to design the optimum transition to MGs for the IEEE 37 and 123 bus test feeders [33]. Candidate power lines to be installed are taken from [34] for the 37, and from [12] for the 123 bus system and the installation cost is 150000 $/mile, discounted over 40 years [35].

Demand shape and solar irradiation profiles are provided by a US Department of Energy database for Crescent City, California [36] and clustered into 8 typical days with hourly time steps using the k-medoids method [24]. Crescent city is located at 41.8°N, 124.2°W on the pacific coast with a yearly global horizontal irradiance of 1440 kWh/m$^2$, and a mean temperature of 10°C [36].

To achieve a similar share of ~40% residential and ~60% commercial demands, all demands above 42 kW peak for the 37 bus, and above 20kW peak for the 123 bus system are defined as commercial.

As PV installations are often customer owned, their capacity is not a decision variable for this case study. PV areas are sampled from a probability distribution relating peak demands and PV capacities, generated from the Pecan Street dataset [37], resulting in 7.8% of the yearly energy demand being supplied by PV, representing the Californian share of 7.6% in 2017 [38]. Panel azimuth and tilt angles are randomly sampled from a list of real PV installations in California [39].

Electricity prices of 15 ¢/kWh for import, 7 ¢/kWh for export are taken from [40] and the reactive power charge of 0.06 ¢/kvarh is based on a rate for customers larger than 500 kW in Crescent City [41].

Techno-economic data of Li-ion storage ($DOD_{max} = 0.85, \eta_{self} = 0.99, \eta_{ch} = \eta_d = 0.98, f_C^{Str} = 1/3, n_{Cyc,max} = 1/day$ [18]) and biogas-fueled micro-turbine DG ($f_{P,max}^{DG} = f_{Q,max}^{DG} = 1, PF_{min}^{DG} = 0$) are taken from institutional reports and shown in Table 1.

*Table 1. Storage and DG parameters*

|  | Storage [18] | DG [42], [43] |
|---|---|---|
| Fixed Cost $c_{Inv,0}^{DER}$ ($) | 87360 | 70250 |
| Var. Cost $c_{Inv,S}^{DER}$ ($/kW) | 670 | 2430 |
| O&M Cost $c_P^{DER}$ ($/kWh) | 0 | 0.122 |

| Q cost $c_Q^{DER}$ (¢/kvarh) | 0.04 [27] | 0.04 [27] |
| Life time $t_{Lf}$ (y) | 15 | 13.3 |

To identify islanding events, [8] suggests the use of Major Event Days (MED), where SAIDI exceeds a threshold identified using the 2.5 beta method. In Crescent city for the time frame 2006 – 2016, the average number of MED is 2 per year [44], which results in $p_i = \lambda_i \Delta t = 2.283 * 10^{-4}$ [45], where $\Delta t$ is the time step of 1 hour. The probabilities of islanding event durations $p_{t_i}$ are computed from the empirical distribution of major event durations between 2006 and 2017 for three Californian utilities [44], [46], [47]. To this aim, the generalized extreme value (GEV) distribution is fitted to the empirical data with a p-value of 0.09 in a chi-squared goodness of fit test with 5% significance level, and used to calculate $p_{t_i}$. Although event durations up to 50 hours occur, MGs are initially designed to island for 24 hours [6] by scaling $p_{t_i}$ to sum to 1 after 24 hours. A sensitivity analysis on islanding duration is performed.

Component reliability data are reported in Table 2. Energy not supplied is penalized with 370 and 3.3 $/kWh for commercial and residential demand, respectively [48]. To account for additional communication, protection and control equipment when transitioning to MGs, islanding investment costs are set to 2 $/MWh demand [3].

*Table 2. Reliability data*

|  | $\lambda$ (1/y) [49] | $\tau$ (h) [50] |
| --- | --- | --- |
| Cable | 0.1 (1/y/mile) | 4 |
| Cable splices | 0.03 | 4 |

The case study optimizations are run on a desktop computer with a 3.6 GHz Intel® Core™ i7-4790 CPU and 16 GB RAM. The optimization problem (1) is modelled with YALMIP [51] in MATLAB R2018a [52] and solved using Gurobi 8.0 [53] with 0.5% optimality gap.

## 4. Results and Discussion

To demonstrate the effectiveness of the proposed approach and to study the influence of resilience and reliability on optimal MG design, the following four alternative design cases are compared for the IEEE 37 bus system shown in Fig. 4, incrementally building the full optimization:

1. minimize investment and operations costs as a base case for the current system ("Base", Levels 1 and 2 in Fig. 1), ignoring reliability and resilience;

2. add reliability costs to Case 1 ("Reliability", Level 4);

3. add islanding costs to Case 1 ("Resilience", Level 3);

4. solve the full problem (1), which is the proposed approach ("Full", Levels 1 – 4).

For Case 1 – "Base" and Case 2 – "Reliability" where islanding is not considered, the impact of islanding events in terms of cost and EENS is calculated ex post, assuming 100% demand not supplied during the events. SAIFI and SAIDI of all solutions are calculated using the method in [28].

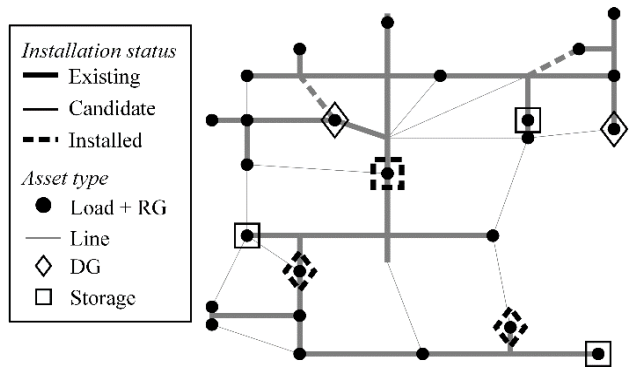

*Fig. 4. IEEE 37 bus topology, peak demands and RG feed-in. Candidate and installed technologies for Case 4 – "Full".*

## 4.1. Design choices and cost

Optimization results are compared in terms of objective values and installed capacities, respectively, in Fig. 5 and Fig. 6. Considering reliability (Case 2) leads to investing in redundant lines, but no DER. Case 3 – "Resilience" enables the MG to island by installing DG and storage, significantly reducing the cost of energy not supplied. DG operating costs are lower than grid imports, therefore, DG complements solar power in providing most of the energy. Case 4 – "Full" entails the minimum total costs. Compared to Case 3, the reliability cost are reduced by installing two additional lines as shown in Fig. 6, and distributing the DG capacity between two buses (see Fig. 4). Comparing the DG and storage capacities in Case 4 to the PCC capacity in Case 1 shows that installing DG capacity equal to 83% of the peak demand complemented with storage rated at 17% of the peak demand, is sufficient to enable islanding for 24 hours.

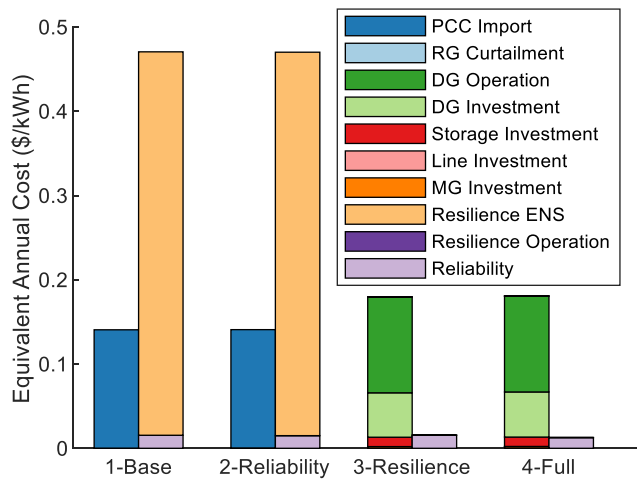

*Fig. 5. Equivalent annual costs per kWh of MG load for Case 1-4, decomposed into investment & operation contribution (left stacks), and resilience and reliability-related contribution (right stacks).*

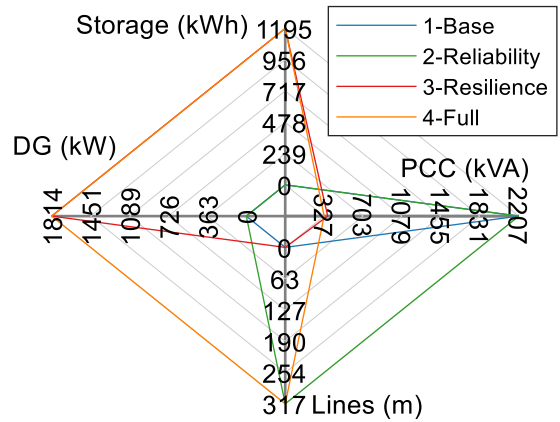

*Fig. 6. Installed lines, storage and DG capacity, and peak power at the PCC for Case 1 – 4.*

### *4.2. Reliability and resilience*

The system reliability and resilience for the three cases is shown in Fig. 7 in terms of SAIFI, SAIDI and EENS. For Case 1 and 2, islanding events are the main contributors to SAIFI, SAIDI and EENS. In Crescent City, the measured SAIFI and SAIDI are 2.57 outages per year and 9.03 hours per year, respectively, averaged in the range 2011 – 2016 [44]. For comparison and validation of the developed methodology, the SAIFI computed for Case 1 shows an excellent agreement with these statistics, whereas the SAIDI is overestimated by 34%. The discrepancies originate from taking islanding event durations from major events across California, due to a lack of data for Crescent City. Despite this inaccuracy, significant reliability and resilience improvements of installing a MG in Case 3 and 4 can be achieved, both compared to the reported data and to Case 1. Case 4 has a slightly better reliability performance than Case 3 due to the redundant lines and the DG placed at two buses (located in the lower part of Fig. 4), which enable self-supply. Compared to the base case, SAIFI, SAIDI and EENS are reduced by 94 to 96% in Case 4.

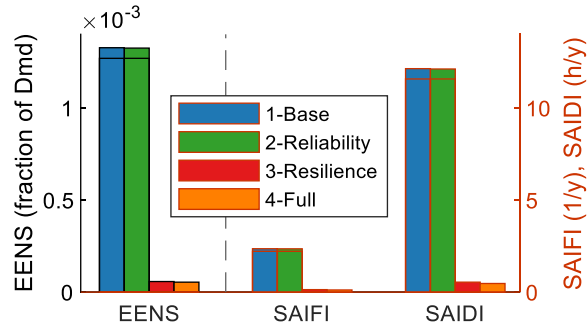

*Fig. 7. Reliability indicators with contributions from islanding (lower) and faults within the MG (upper stack elements).*

### 4.3. Sensitivity Analysis

To generalize the presented case study for other design conditions/parameters and to measure the robustness of the identified optimal design choices, sensitivity analyses are conducted with respect to the following parameters in the full problem (1):

1. the islanding duration is varied in 6 hour increments between 6 and 30 hours to capture different islanding events [7], [8];

2. the installed renewable generator capacity varies from 50% to 200% of the nominal case, representing current feeders with low PV penetration as well as a high PV scenario in 2030, where the US renewable generation is expected to double compared to 2017 [54];

3. storage investment costs are reduced in 10% increments by up to 40%, reflecting possible battery technology advances until 2030 [55].

Reducing storage investment costs has no effect on the optimal storage and DG capacity, with the exception of the 6 hour islanding at 50% RG capacity conditions, which result in the installation of 2.6 MWh of storage capacity. This invariance indicates that the main driver for DER investments is to provide generation capacity for islanding conditions and that the relative cost of storage and DG does not influence this balance. Therefore, Fig. 8 reports the installed DG and storage capacity for variable islanding duration and RG capacity and for the nominal value of the storage cost reported in Table 1.

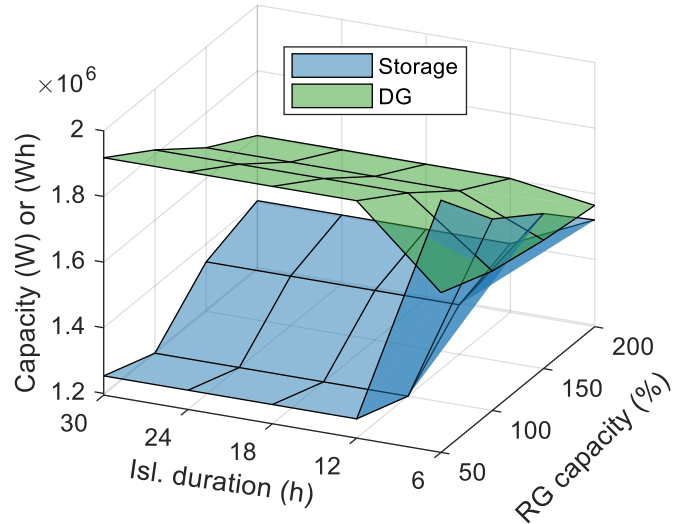

*Fig. 8. Sensitivity analysis of storage and DG capacity to maximum islanding event duration and RG capacity.*

With increasing RG capacity, storage replaces DG capacity because the MG can buffer and use more RG power during islanding situations. The design invariance for islanding durations over 12 hours demonstrates that if a MG is able to island for 12 hours, it can also island for 30 hours, because sufficient DER are installed to supply all demands. If the MG is designed for 6 hour islanding events, however, the DG capacity is replaced by storage, because storage can be effectively used to provide energy for such short durations. MGs designed for a maximum islanding duration of 6 hours are unable to satisfy all demands during longer interruptions, showing the importance of designing MGs for events longer than the proposed durations in [11], [14], if resilience is required.

### *4.4. Multi-Microgrids*

The IEEE 123 bus case is investigated to determine the correlation between system size and computational efforts, and to explore the cost and reliability performance of dividing one feeder into multiple MGs. Starting from a single MG, the IEEE 123 bus feeder is divided into up to 4 MGs, by introducing MG partitions A, B, and C in Fig. 9.

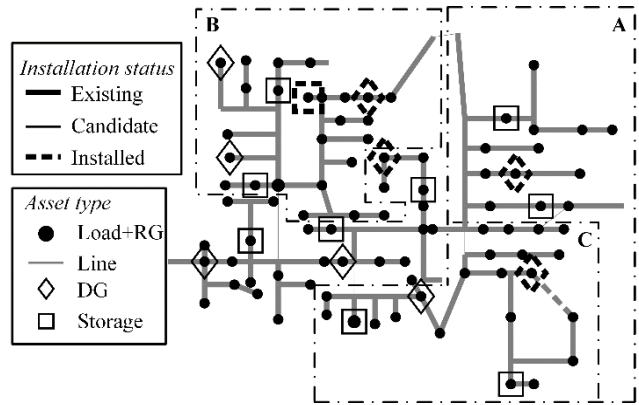

*Fig. 9. IEEE 123 bus topology, partitions for the multi-MG cases and installed technologies for the four MG case. The four MGs partitions are identified by the dash-dotted contours.*

The resulting MGs and optimization problem sizes are shown in Table 3. The single MG system is made of MG I, the 2 MG system of II and III, the 3 MG system of II, IV and V, and the 4 MG system of IV, V, VI and VII.

*Table 3. MG optimization problem sizes resulting from partitioning the IEEE 123 bus system.*

| MG | Buses | Continuous variables ($10^6$) | Binary variables ($10^3$) | Constraints ($10^6$) |
|---|---|---|---|---|
| I Full system | 123 | 8.59 | 21.8 | 51.2 |
| II Partition A | 51 | 3.06 | 3.82 | 19.3 |
| III All buses except A | 72 | 4.42 | 7.07 | 27.1 |
| IV Partition B | 35 | 1.99 | 1.69 | 12.5 |
| V All buses except A & B | 37 | 2.06 | 1.81 | 12.9 |
| VI Partition C | 27 | 1.44 | 0.92 | 9.29 |
| VII Partition A without C | 24 | 1.41 | 0.95 | 8.97 |

Fig. 10 shows the total equivalent annual cost from the 1MG to the 4MG case, including the original, passive grid without DER ("No MG"). The most significant change happens between No MG and 1 MG, where the investment and operation cost increase by 25%, while the reliability and resilience in terms of EENS is reduced by 96%. In the 2 – 4 MG cases, horizontal lines within the bar stacks show the costs are split between MGs. With increasing numbers of MGs, costs are shifted from interruption (purple) to investment and operation because smaller parts of the feeder can island separately during upstream faults, but each MG requires separate DER. Furthermore, storage (red) is substituted by DG (light green) with increasing numbers of MGs. In the 4 MG case, only MG IV installs storage, resulting in the smallest overall cost.

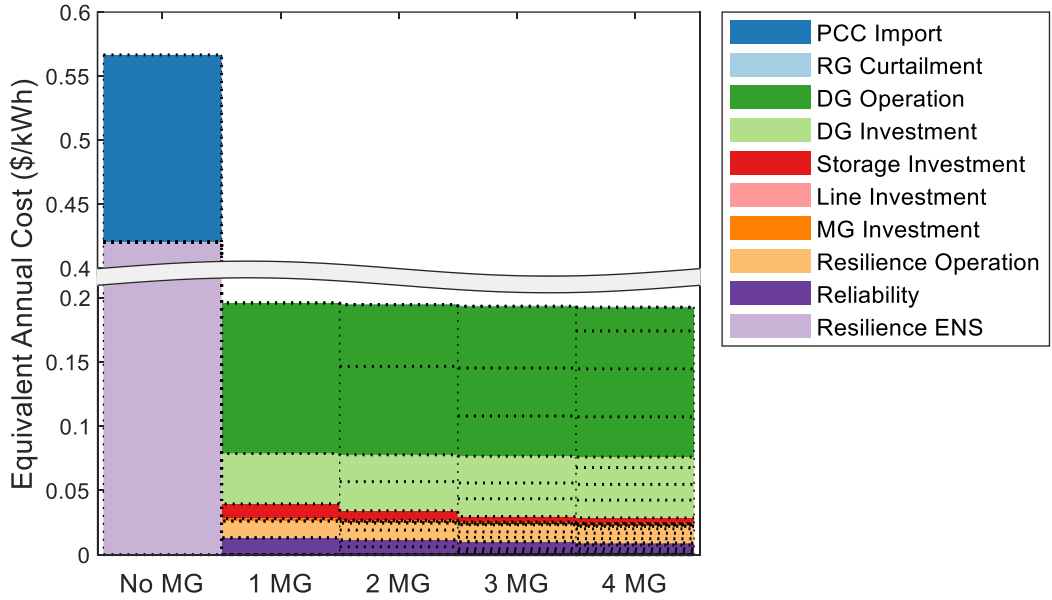

*Fig. 10. Equivalent Annual Cost of the 123 bus feeder, split into up to 4 MGs including the passive system ("No MG") for comparison. The cost for individual MGs are shown by horizontal dotted lines, which split each cost for each MG partition.*

Fig. 11 shows the computational time required to solve the full design optimization (1) for different MG sizes. Circle markers indicate CCG solver times on the desktop computer from Section 3. As a reference, the same optimizations are solved on the Euler cluster [56] using Benders decomposition [57] in CPLEX 12.8 [58] with 8 cores and up to 137 GB RAM. CCG outperforms Benders decomposition by up to a factor of 12, solving the single MG case in 5 hours and proving that the proposed approach is scalable for larger systems in contrast to Benders decomposition.

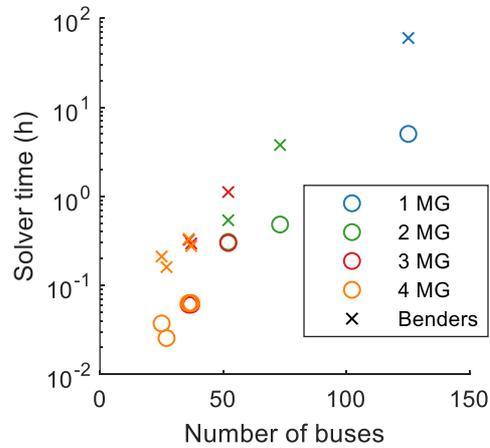

*Fig. 11. Solver times for the 123 bus case, split into up to 4 MGs. Comparison of CCG (circles) and Benders decomposition (x).*

## 5. Conclusion

This paper investigates optimal MG design for techno-economic and reliability objectives, modelling and minimizing costs of energy not supplied due to islanding events and faults within the MG. The effectiveness of the proposed method is demonstrated on modified IEEE 37 and 123 bus test feeders, where investment, operation and interruption costs can be minimized by installing and optimally operating storage and DG.

This study shows that combining design, grid-tied operations, islanding and faults yields an optimal transition to highly reliable and resilient MGs with minimum total costs.

In particular, the incorporation of stochastic islanding events into the operational optimization allows for effective anticipation of such events. Additionally, the proposed CCG method ensures the scalability of the optimization framework to large design problems.

The reliability and resilience improvements of transitioning distribution grids to MGs are demonstrated with reference to two case studies. Storage is cost-optimal in all cases, particularly with high RG penetration, and should, therefore, be included in MG design optimization.

Considering only reliability or resilience can lead to higher demand not supplied, resulting in increased social cost. Designing MGs for 12 hours of islanded operation can be sufficient to ensure islanding capabilities for at least 30 hours, which highlights the benefits and importance of considering longer islanding durations to design resilient MGs.

Finally, the proposed islanding and grid-tied operation model can be adapted for optimal scheduling of resilient MGs by introducing forecast uncertainty of demands and RG as well as demand side management, which is considered to be future work.